\documentclass[12pt,reqno]{amsart}
\usepackage{amssymb}

\input xy
\xyoption{all}

\renewcommand{\epsilon}{\varepsilon}

\renewcommand{\emptyset}{\varnothing}
\renewcommand{\phi}{\varphi}

\title[spaces with  stabilizers of prime power order]{On groups acting
  on contractible spaces with stabilizers of prime power order}

\author{Ian J. Leary}
\address{Department of Mathematics, The Ohio State University, 
231 W 18th Ave, Columbus Ohio 43210}
\address{Heilbronn Institute, University of Bristol, Royal Fort
  Annexe, Bristol BS8 1TW}
\email{leary@math.ohio-state.edu}
\author{Brita E. A. Nucinkis}
\address{School of Mathematics, University of Southampton, Southampton
  SO17 1BJ, UK} 
\email{B.E.A.Nucinkis@soton.ac.uk}
\thanks{The first named author was partially supported by NSF grant
  DMS-0804226 and by the Heilbronn Institute.}

\newtheorem{theorem}{Theorem}[section]
\newtheorem{proposition}[theorem]{Proposition}
\newtheorem{corollary}[theorem]{Corollary}

\newtheorem{conjecture}[theorem]{Conjecture}

\theoremstyle{definition}

\theoremstyle{remark}
\newtheorem{remark}[theorem]{Remark}

\newcommand{\clh}{{\scriptstyle\bf H}}

\newcommand{\FF}{\mathfrak F}
\newcommand{\ff}{\mathcal F}

\newcommand{\PP}{\mathfrak P}
\newcommand{\OO}{\mathfrak O}

\newcommand{\Z}{\mathbb Z}

\newcommand{\F}{\mathbb F}

\newcommand{\FP}{\mathrm{FP}}

\newcommand{\FFP}{{\mathcal F}{\operatorname{FP}}}

\newcommand{\E}{\mathrm{E}}

\newcommand{\cohom}[3]{H^{{\raise1pt\hbox{$\scriptstyle#1$}}}(#2\>\!,#3)}
\newcommand{\tatecohom}[3]{\widehat H^{{\raise1pt\hbox{$\scriptstyle#1$}}}(#2\>\!,#3)}

\newcommand{\Cohom}[3]{H^{{\raise1pt\hbox{$\scriptstyle#1$}}}\big(#2\>\!,#3\big)}
\newcommand{\Tatecohom}[3]{\widehat H^{{\raise1pt\hbox{$\scriptstyle#1$}}}\big(#2\>\!,#3\big)}

\newcommand{\homol}[3]{H_{{\lower1pt\hbox{$\scriptstyle#1$}}}(#2\>\!,#3)}
\newcommand{\homolog}[2]{H_{{\lower1pt\hbox{$\scriptstyle#1$}}}(#2)}

\newcommand{\IND}{\operatorname{Ind}}

\newcommand{\ind}[3]{{\IND}_{#1}^{#2}#3}

\newcommand{\mono}{\rightarrowtail}
\newcommand{\epi}{\twoheadrightarrow}

\newcommand{\Epi}[1]{\buildrel {#1}\over\twoheadrightarrow}

\newcommand{\EG}{{{\mathrm{E}}G}}
\newcommand{\eg}{{\underline{\mathrm{E}}G}}
\newcommand{\efg}{{{\mathrm{E}}_{\FF}G}}
\newcommand{\epg}{{\mathrm{E}}_{\PP}G}
\newcommand{\evcg}{{\operatorname E}_{\mathcal{VC}}G}

\newcommand{\singp}{{\mathrm{sing}(p)}}
\newcommand{\Hom}{\operatorname{Hom}}

\newcommand{\D}{\Delta}

\newcommand{\hf}{{\clh \mathfrak F}}
\newcommand{\hp}{{\clh\mathfrak P}}

\newcommand{\HP}[1]{{\clh}_{#1}{\mathfrak P}}
\newcommand{\HF}[1]{{\clh}_{#1}{\mathfrak F}}

\begin{document}

\begin{abstract}
Let $\FF$ denote the class of finite groups, and let $\PP$ denote the
subclass consisting of groups of prime power order.  We study group
actions on topological spaces in which either (1) all stabilizers lie
in $\PP$ or (2) all stabilizers lie in $\FF$.  We compare the
classifying spaces for actions with stabilizers in $\FF$ and $\PP$,
the Kropholler hierarchies built on $\FF$ and $\PP$, and group
cohomology relative to $\FF$ and to $\PP$.  In terms of standard
notations, we show that $\FF \subset \HP1\subset \HF1$, with all 
inclusions proper; that $\hf=\hp$; that $\FF H^*(G;-)=\PP H^*(G;-)$; 
and that $\epg$ is finite-dimensional if and only if $\efg$ is 
finite-dimensional and every finite subgroup of $G$ is in $\PP$.  
\end{abstract}

\maketitle
\section{Introduction}\label{intro}

Let $\ff$ denote a class of groups, by which we
mean a collection of groups which is closed under isomorphism and taking subgroups.
A $G$-CW-complex $X$ is said to be a model for ${{\E}_\ff}G$, 
the classifying space for actions of $G$ with stabilizers in $\ff$, 
if the fixed point set $X^H$ is contractible for $H\in \ff$ and is 
empty for $H\notin \ff$.  The most common classes considered are 
the class of trivial groups and the class $\FF$ 
consisting of all finite groups.  In these cases ${{\E}_\ff}G$ 
is often denoted ${\E}G$ and $\eg$ respectively.  Note that $\EG$ is the 
total space of the universal principal $G$-bundle, or equivalently 
the universal covering space of an Eilenberg-Mac~Lane space for~$G$.  
The space $\eg$ is called the classifying space for proper actions 
of~$G$.  Recently there has been much interest in finiteness
conditions for the spaces ${{\E}_\ff}G$, especially for $\eg$.  
Milnor and Segal's constructions of $\EG$ both generalize easily to 
construct models for any ${{\E}_\ff}G$, and one can show that any two 
models for ${{\E}_\ff}G$ are naturally equivariantly homotopy equivalent.

For some purposes the structure of the fixed point sets for subgroups
in $\ff$ is irrelevant.  For example, a group is in Kropholler's class 
$\clh_1\ff$ if there is any finite-dimensional contractible
$G$-CW-complex $X$ with all stabilizers in $\ff$.  The class $\clh_1
\ff$ is the first stage of a hierarchy whose union is Kropholler's 
class $\clh \ff$ of hierarchically decomposable groups~\cite{k1}.  
(These definitions were first considered for the class $\FF$ of 
all finite groups, but work for any class $\ff$.)  

A priori, the class $\clh_1\ff$ may contain groups $G$ that do not admit a 
finite-dimensional model for ${{\E}_\ff}G$, and we shall give such
examples in the case when $\ff=\PP$, the class of groups of prime
power order.  By contrast, in the case when $\ff = \FF$, no group 
$G$ is known to lie in $\clh_1\FF$ without also admitting a
finite-dimensional model for $\eg$.  A construction due to Serre shows
that every group $G$ in $\clh_1\FF$ that is virtually torsion-free has
a finite-dimensional $\eg$~\cite{brbo}, and the authors have given
examples of $G$ for which the minimal dimension of a contractible
$G$-CW-complex is lower than the minimal dimension of a model for 
$\eg$~\cite{vfg}.  These examples $G$ also have the property 
that they admit a contractible $G$-CW-complex with finitely many 
orbits of cells, but that they do not admit any model for $\eg$ 
with finitely many orbits of cells.  

Throughout this paper, $\FF$ will denote the class of finite groups, 
and $\PP$ will denote the class of finite groups of prime power
order.  We compare the classifying space for $G$-actions with stabilizers 
in $\PP$ with the more well-known $\eg$, and we compare the Kropholler 
hierarchies built on $\FF$ and $\PP$.  We show that a finite group 
$G$ that is not of prime power order cannot admit a
finite-dimensional $\epg$, but that every finite group is in 
$\clh_1\PP$.  We also construct a group that is in $\clh_1\FF$ but not 
in $\clh_1\PP$, and we show that $\clh\PP=\clh\FF$.  

In the final section we shall contrast this with cohomology relative to 
all finite subgroups. The relative cohomological dimension can
be viewed as a generalisation of the virtual cohomological dimension, since
for virtually torsion free groups these are equal, see~\cite{mn}.
By a result of Bouc \cite{bouc, kw} it follows that groups belonging to 
$\HF1$ have finite relative cohomological dimension, but the converse it
not known.  In contrast to our results concerning classifying spaces, 
we show that cohomology relative to subgroups in $\FF$
is naturally isomorphic to cohomology relative to subgroups in $\PP$.

\medskip
\section{Classifying spaces for actions with stabilizers in $\PP$}
\medskip

\begin{theorem}\label{noefg}
Let $G$ be a finite group. Then $G$ has a finite dimensional model for
$\epg$ if and only if $G$ has prime power order.  
\end{theorem}

\noindent{\sl Proof.} If $G$ has prime power order, then a single
point may be taken as a model for $\epg$.  
Now let $G$ be an arbitrary finite group, let $p$ be a prime dividing 
the order of $G$, and 
assume that there is a $p$-subgroup
$P< G$, such that $N_G(P)$ is not a $p$-group. Then the Weyl-group
$WP = N_G(P)/P$ contains a subgroup $H$ of order prime to $p$.
Assume $G$ has a finite dimensional model for $\epg$, $X$ say.
Then the augmented cellular chain complex of the $P$-fixed point set,
$X^P$, is a finite length
resolution of $\Z$ by free $H$-modules.  This gives a contradiction, 
since $\Z$ has infinite projective dimension as an $H$-module for 
any non-trivial finite group $H$.

Therefore we may suppose  that $G$ is not in $\PP$ and 
for all non-trivial subgroups $P\in \PP$, the normalizer $N_G(P)$ is 
also in $\PP$.  Now let $N$ be a minimal normal subgroup
of $G$. This cannot lie in $\PP$ and hence has the same properties as $G$. 
Thus, by minimality we may assume $N =G$ and $G$ is simple. 

Choose two distinct Sylow $p$-subgroups $P$ and $Q$ of $G$, such that
their intersection, $I$ say, is of maximal order, and assume that $I$ is non-trivial.  
Now, the normalizers $N_P(I)$ and $N_Q(I)$ contain $I$ as a proper subgroup.
Also, the group $\langle N_P(I), N_Q(I) \rangle$ does not contain $P$
and $Q$ and neither does $N_G(I) \geq \langle N_P(I), N_Q(I) \rangle$,
which is a $p$-group by assumption.  Hence there exists a
Sylow $p$-subgroup $R$ containing $N_G(I)$ and $R \cap P \geq N_P(I)$.
Thus $|R \cap P| > |I| = |P \cap Q|$, which contradicts the maximality
of $I$.  Therefore we may assume that in $G$ all Sylow $p$-subgroups
intersect trivially.  In such a group we have, for $P$ a
Sylow $p$-subgroup:
$$ H^*(G,\F_p)  \cong H^*(P,\F_p),$$
see for example \cite[Theorem III.10.3]{brbo}.  


Any non-trivial $p$-group has non-trivial abelianization, and hence
$H^1(P,\F_p)$, which is naturally isomorphic to $\Hom(P,\F_p)$, 
is non-trivial.  But this implies 
that $H^1(G,\F_p)\cong \Hom(G,\F_p)$
is non-trivial, and so $G$ admits a surjective homomorphism to 
a group of order $p$.  Since $G$ is not in $\PP$, it follows that 
$G$ cannot be simple, which gives the final contradiction.  
\qed

\begin{remark}
Yoav Segev pointed out that the above proof could be shortened by quoting the 
Frobenius normal $p$-complement theorem~\cite[5.26]{isaacs}.   
\end{remark}

\begin{corollary}
 For a group $G$, the following are equivalent.  
\begin{enumerate} 
\item 
$G$ admits a finite-dimensional $\epg$; 
\item
Every finite subgroup of $G$ is in $\PP$ and $G$ admits a
 finite-dimensional $\eg$. \qed 
\end{enumerate}
 \end{corollary}

\begin{remark} We conclude the section with a remark on the type
 of $\epg$. It can proved analogously to L\"uck's proof for $\eg$
 \cite{lupa} that a group $G$ admits a finite type model for $\epg$ if
 and only if $G$ has finitely many conjugacy classes of groups of
 prime power order and the Weyl-groups $N_G(P)/P$ for all subgroups
 $P$ of prime power order are finitely presented and of type
 $\FP_\infty.$ Hence any group admitting a finite type $\eg$ also
 admits a finite type $\epg$. Recall that a finite extension of a
 group admitting a finite model for $\EG$ always has finitely many
 conjugacy classes of subgroups of prime power order
 \cite[IX.13.2]{brbo}.  Hence the groups exhibited in \cite[Example
 7.4]{vfg} are groups admitting a finite type $\epg$ which do not
 admit a finite type $\eg.$  

This behaviour is in stark contrast to that of $\evcg$, the
classifying space with virtually cyclic isotropy.  Any group admitting
a finite dimensional model for $\evcg$ admits a finite dimensional
model for $\eg$, see \cite{lueckweiermann} and the converse also holds
for a large class of groups including all polycyclic-by-finite and all
hyperbolic groups \cite{jpl, lueckweiermann}. Furthermore, any group
admitting a finite type model for $\evcg$ also admits a finite type
model for $\eg$ \cite{komn}, but it is conjectured \cite{jpl} that any
group admitting a finite model for $\evcg$ has to be virtually
cyclic. This has been shown for a class of groups containing all
hyperbolic groups \cite{jpl} and for elementary amenable groups
\cite{komn}.

\end{remark}

\medskip
\section{The hierarchies $\clh\FF$ and $\clh\PP$}
\medskip

\begin{proposition}\label{prop:hierarchy}
Let $X$ be a finite dimensional contractible $G$-CW-complex such that
all stabilizers are finite.  If there is a bound on the orders of the
stabilizers then there exists a finite dimensional contractible
$G$-CW-complex $Y$ and an equivariant map $f: Y \to X$ such that $Y^H
= \emptyset$ if $H$ is not a $p$-group.
\end{proposition}

\noindent{\sl Proof.}  Using the equivariant form of the simplicial 
approximation theorem, we may assume that $X$ is a simplicial 
$G$-CW-complex.  To simplify notation the phrase `$G$-space' shall 
mean `simplicial $G$-CW-complex' and `$G$-map' will mean
`$G$-equivariant simplicial map' throughout the rest of the proof.  
The space $Y$ will be a $G$-space in this sense and the map
$f:Y\rightarrow X$ will be a $G$-map in this sense.  The $G$-space 
$Y$ is constructed in two stages.  Firstly, for each finite $H\leq G$
we build a finite-dimensional contractible $H$-space $Y_H$ with 
the property that all simplex stabilizers in $Y_H$ lie in $\PP$.  

Suppose for now that each such $H$-space $Y_H$ has been constructed.  
Using the $G$-equivariant
form of the construction used in \cite[Section 8]{km} the space $Y$ is
constructed as follows.  Let $I$ be an indexing set for the $G$-orbits of 
vertices in $X$.  For each $i\in I$, let $v_i$ be a representative of
the corresponding orbit, and let $H_i$ be the stabilizer of $v_i$.  
Let $X^0$ denote the 0-skeleton of $X$.  Define a $G$-space $Y^0$ by 
$$Y^0 = \coprod_{i\in I} G\times_{H_i} Y_{H_i},$$ 
and define a $G$-map $f:Y^0\rightarrow X^0$ by $f(g,y)= g.v_i$ for all 
$i\in I$, for all $g\in G$ and for all $y\in Y_{H_i}$.  For each
vertex $w$ of $X$, let $Y(w)= f^{-1}(w)\subset Y^0$.  Each $Y(w)$ is 
a contractible subspace of $Y^0$, and the stabilizer of $w$ acts on 
$Y(w)$.  

Now for $\sigma= (w_0,\ldots,w_n)$ an $n$-simplex of $X$, define a 
space $Y(\sigma)$ as the join 
$$Y(\sigma)= Y(w_0)*Y(w_1)*\cdots *Y(w_n).$$ Each vertex of
$Y(\sigma)$ is already a vertex of one of the $Y(w_i)$, and so the map
$f:Y^0\rightarrow X^0$ defines a unique simplicial map
$f:Y(\sigma)\rightarrow \sigma$.  By construction, whenever $\tau$ is
a face of $\sigma$, the space $Y(\tau)$ is identified with a subspace
of $Y(\sigma)$.  This allows us to define $Y$ and $f:Y\rightarrow X$
as the colimit over the simplices $\sigma$ of $X$ of the subspaces
$Y(\sigma)$, and to define $f:Y\rightarrow X$, which is a $G$-map of
$G$-spaces.  Since each $Y(\sigma)$ is contractible, it follows that
$f$ is a homotopy equivalence, and hence $Y$ is also contractible (see
\cite[Corollary 8.6]{km}).

It remains to build the $H$-space $Y_H$ for each finite group $H<G$.
In the case when $H\in \PP$ we may take a single point to be $Y_H$,
and so we may suppose that $H\notin\PP$.  Fix such a subgroup $H$, and
suppose that we are able to construct a finite-dimensional
contractible $H$-space $Z_H$ in which each stabilizer is a proper
subgroup of $H$.  We may assume by induction that for each $K<H$ we
have already constructed the $K$-space $Y_K$.  The $H$-space $Y_H$ can
now be constructed from $Z_H$ and the spaces $Y_K$ using a process
similar to the construction of $Y$ from $X$ and the spaces $Y_H$. 
It remains to construct the $H$-space $Z_H$.  

An explicit construction of an $H$-space $Z_H$ with the required
properties is given in \cite{ivf}.  We therefore provide only a sketch
of the argument.  We may assume that $H$ is not in $\PP$.  Let $S$ be
the unit sphere in the reduced regular complex representation of $H$,
so that $S$ is a topological space with $H$-action such that the
stabilizer of every point of $S$ is a proper subgroup of $H$.  Since
$H$ is not in $\PP$, there are $H$-orbits in $S$ of coprime lengths.
Using this property, it can be shown that the sphere $S$ admits an
$H$-equivariant self-map $g:S\rightarrow S$ of degree zero.  The
$H$-space $Z_H$ is defined to be the infinite mapping telescope
(suitably triangulated) of the map $g$.  \qed

\begin{corollary}\label{cor:hier}
If $G$ is in $\clh_1\FF$ and there is a bound on the orders of the 
finite subgroups of $G$, then $G$ is in $\clh_1\PP$.  \qed
\end{corollary}

\begin{remark} In Proposition~\ref{prop:hierarchy}, 
the bound on the orders of the stabilizers of $X$ is used only to 
give a bound on the dimensions of the spaces $Y_H$.  In 
Theorem~\ref{cor:hierarchy} we shall show that $\clh_1\FF\neq 
\clh_1\PP$.  
\end{remark}


\begin{remark} The construction in Proposition~\ref{prop:hierarchy} 
does not preserve cocompactness, because for most finite groups $H$,
the space $Y_H$ used in the construction cannot be chosen to be
finite.  A result similar to Proposition~\ref{prop:hierarchy} but
preserving cocompactness can be obtained by replacing $\PP$ by a
larger class $\OO$ of groups.  Here $\OO$ is defined to be the class
of $\PP$-by-cyclic-by-$\PP$-groups.  A theorem of Oliver~\cite{oliver}
implies that any finite group $H$ that is not in $\OO$ admits a
\emph{finite} contractible $H$-CW-complex $Z'_H$ in which all
stabilizers are proper subgroups of $H$.  Applying the same argument
as in the proof of Proposition~\ref{prop:hierarchy}, one can show that
given any contractible $G$-CW-complex $X$ with all stabilizers in
$\FF$, there is a contractible $G$-CW-complex $Y'$ with all
stabilizers in $\OO$ and a proper equivariant map $f':Y'\rightarrow
X$.  (By proper, we mean that the inverse image of any compact subset
of $X$ is compact.)
\end{remark}


For $X$ a $G$-CW-complex with stabilizers in $\FF$, and $p$ a prime, let 
$X_\singp$ denote the subcomplex consisting of points whose stabilizer 
has order divisible by $p$.  For $G$ a group and $p$ a prime, let
$S_p(G)$ denote the poset of non-trivial finite $p$-subgroups of $G$.  

\begin{proposition} \label{prop:sing}
Suppose that $X$ is a finite-dimensional contractible 
$G$-CW-complex with all stabilizers in $\PP$.  For each prime $p$, 
the mod-$p$ homology of $X_\singp$ is isomorphic to the mod-$p$ 
homology of the (realization of the) poset $S_p(G)$.  
\end{proposition}


\noindent{\sl Proof.} Fix a prime $p$, and to simplify notation let 
$S$ denote the realization of the poset $S_p(G)$.  For $P$ a
non-trivial $p$-subgroup of $G$, let $X^P$ denote the points fixed by
$P$, and let $S_{\geq P}$ denote the realization of the subposet of 
$S_p(G)$ consisting of all $p$-subgroups that contain $P$.  By the 
P. A. Smith theorem~\cite{bred}, each $X^P$ is mod-$p$ acyclic.  
Each $S_{\geq P}$ is contractible since it is equal to a cone with 
apex $P$.  Let $P$ and $Q$ be $p$-subgroups of $G$, and let $R=\langle
P, Q\rangle$, the subgroup of $G$ generated by $P$ and $Q$.  If $R$ is
a $p$-group then $X^P\cap X^Q=X^R$, and otherwise $X^P\cap X^Q$ is
empty.  Similarly, $S_{\geq P}\cap S_{\geq Q}= S_{\geq R}$ if $R$ is a 
$p$-group and $S_{\geq P}\cap S_{\geq Q}$ is empty if $R$ is not a
$p$-group.  

Since each $X^P$ is mod-$p$ acyclic, the mod-$p$ homology
$H_*(X_\singp)$ is isomorphic to the mod-$p$ homology of the nerve of
the covering $X_\singp=\bigcup_P X^P$.  Similarly, the mod-$p$ homology 
$H_*(S_p(G))$ is isomorphic to the mod-$p$ homology of the nerve of 
the covering $S_p(G)=\bigcup_P S_{\geq P}$.  By the remarks in the 
first paragraph, these two nerves are isomorphic.  \qed 


\begin{proposition} \label{prop:lowerbd}
Let $k$ be a finite field, and let $G$ be the group of $k$-points of a 
reductive algebraic group over $k$, whose commutator subgroup has $k$-rank $n$.  
(For example, $G= SL_{n+1}(k)$, or $GL_{n+1}(k)$.)  Any finite-dimensional 
contractible $G$-CW-complex 
with stabilizers in $\PP$ has dimension at least $n$.
\end{proposition} 

\noindent {\sl Proof.} The hypotheses on $G$ imply that $G$ acts on a
spherical building $\Delta$ of dimension $n-1$~\cite[Appendix on algebraic
groups]{abramenkobrown,brown}.  Any such building is homotopy
equivalent to a wedge of $(n-1)$-spheres.  
Quillen has shown that $\Delta$ 
is homotopy equivalent to the realization of $S_p(G)$, where $p$ is
the characteristic of the field $k$~\cite[Proposition~2.1 and 
  Theorem~3.1]{quillen}.  It follows that $S_p(G)$ is homotopy
equivalent to a wedge of $(n-1)$-spheres, and in particular 
the mod-$p$ homology group $H_{n-1}(S_p(G))$ is non-zero.  

Now suppose that $X$ is a finite-dimensional 
contractible $G$-CW-complex with stabilizers in $\PP$.  
Using Proposition~\ref{prop:sing}, one sees that the mod-$p$ 
homology group $H_{n-1}(X_\singp)$ is non-zero.  It 
follows that $X$ must have dimension at least $n$.  
\qed 

\begin{remark} 
In \cite{segev} it is shown that if $G$ is a finite simple group of
Lie type, of Lie rank $n$, then any contractible $G$-CW-complex of
dimension strictly less than $n$ contains a point fixed by $G$.
(Theorem~1 of \cite{segev} contains the additional hypothesis that the
$G$-CW-complex should be finite, but this is not used in the proof.)
A similar argument to that used in \cite[Theorem~2]{segev} was used in
\cite{ivf} to show that when $G=SL_{n+1}(\F_p)$, every contractible
$G$-CW-complex without a global fixed point has dimension at least
$n$.  Note that Proposition~\ref{prop:lowerbd} applies in greater
generality than these results.  For example, the Conner-Floyd
construction \cite{cf} shows that whenever the multiplicative group of
$k$ does not have prime-power order, there is, for any $n \geq 1$, a
4-dimensional contractible $Gl_n(k)$-CW-complex without a global fixed
point.
\end{remark}

\begin{theorem}\label{cor:hierarchy}
There are the following strict containments and equalities between
classes of groups: 
\begin{enumerate}
\item $\FF\subsetneq \HP1$; 
\item $\HP1 \subsetneq \HF1$; 
\item $\hf = \hp$.  

\end{enumerate}
\end{theorem}

\noindent{\sl Proof.}  Corollary~\ref{cor:hier} shows that
$\FF\subseteq \HP1$.  The free product of two cyclic groups of prime
order is in $\HP1$ and is not finite.  The claim that $\hf=\hp$
follows from the inequalities $\PP\subseteq \FF\subseteq \HP1$, 
and the claim $\HP1\subseteq \HF1$ follows from $\PP\subseteq \FF$.  

It remains to exhibit a group $G$ that is in $\HF1$ but not in
$\HP1$.  Let $G=SL_\infty(\F_p)$, the direct limit of the groups 
$G_n=SL_n(\F_p)$, where $G_n$ is included in $G_{n+1}$ as the `top 
corner'.  As a countable locally-finite group, $G$ acts with finite 
stabilizers on a tree.  (Explicitly, the vertex set $V$ and edge set 
$E$ are both equal as $G$-sets to the disjoint union of the sets of cosets 
$G/G_1\cup G/G_2\cup \cdots$, with the edge $gG_i$ joining the vertex 
$gG_i$ to the vertex $gG_{i+1}$.)  It follows that $G\in \HF1$.  
By Proposition~\ref{prop:lowerbd}, $G$ cannot be in $\HP1$.  
\qed

%

\begin{remark}
Let $G$ be a group in $\hf$ that is also of type $\FP_\infty$.
By a result of Kropholler~\cite{k1}, there is a bound on the orders 
of finite subgroups of $G$, and Kropholler-Mislin show that $G$ is 
in $\HF1$~\cite{km}.  Corollary~\ref{cor:hier} shows that $G$ is in 
$\HP1$.  
\end{remark}

\medskip
\section{Cohomology relative to a class of groups}
\medskip

Let $\Delta$ denote a $G$-set, and let $\Z\Delta$ denote the
corresponding $G$-module.  For $\delta\in \Delta$, we write $G_\delta$
for the stabilizer of $\delta$.  A short exact sequence $A\mono B\epi
C$ of $G$-modules is said to be {\sl $\Delta$-split} if and only if it
splits as a sequence of $G_\delta$-modules for each $\delta\in\Delta$.
Equivalently, the sequence is $\Delta$-split if and only if the
following sequence of $\Z G$-modules splits: $A\otimes \Z\Delta \mono
B \otimes \Z\Delta \epi C \otimes \Z\Delta$ \cite{N1}.

We say a $G$ module is $\Delta$-projective if it is a direct summand
of a $G$-module of the form $N\otimes \Z\Delta$, where $N$ is an
arbitrary $G$-module. $\Delta$-projectives satisfy analogue properties
to ordinary projectives.  Furthermore, for each $\delta$, and each
$G_\delta$-module $M$, the induced module $\ind{{G}_\delta}{G}M$ is
$\Delta$-projective.  Given two $G$-sets $\Delta_1$ and $\Delta_2$ and
a $G$-map $\Delta_1 \to \Delta_2$ then $\Delta_1$-projectives are
$\Delta_2$-projective and $\Delta_2$-split sequences are
$\Delta_1$-split. For more detail the reader is referred to \cite{N1}. 

Now suppose that $\ff$ is a class of groups closed under taking subgroups.   
We consider $G$-sets $\Delta$ satisfying the following condition, for all $H\leq G$: 
$$\Delta^H \neq \emptyset \iff H \in {\ff} . \leqno(*)$$ 
There are $G$-maps between any two $G$-sets satisfying condition
($\ast$), and so we may define an $\ff$-projective module to be a
$\Delta$-split module for any such $\Delta$.  Similarly, an
$\ff$-split exact sequence of $G$-modules is defined to be a 
$\Delta$-split sequence.  If $\Delta$ satisfies ($\ast$) and $M$ 
is any $G$-module, the module $M\otimes \Z\Delta$ is $\ff$-projective
and admits an $\ff$-split surjection to $M$.  
This leads to a construction of homology relative to $\ff$.    
An $\ff$-projective resolution of a module
$M$ is an $\ff$-split exact sequence 
$$ \cdots \to P_{n+1} \to P_n \to \cdots \to P_0 \to M \to 0,$$
where all $P_i$ are $\ff$-projective.  Group cohomology relative to 
$\ff$, denoted $\ff H^*(G;N)$ can now be defined as the cohomology 
of the cochain complex $\Hom_G(P_*,N)$, where $P_*$ is an
$\ff$-projective resolution of $\Z$.  

We say that a module $M$ is of type 
$\ff\FP_n$ if $M$ admits an $\ff$-projective resolution in which $P_i$
is finitely generated for $0\leq i\leq n$.  It has been shown that 
modules of type $\ff\FP_n$ are of type $\FP_n$ \cite{N1}. We will say 
that a group $G$ is of type $\ff\FP_n$ if the trivial $G$-module $\Z$ 
is of type $\ff\FP_n$.

We now specialize to the cases when $\ff=\FF$ and $\ff=\PP$.  

\begin{theorem} The following properties hold.  
\begin{enumerate}
\item A short exact sequence of $G$-modules is $\FF$-split if and 
only if it is $\PP$-split.
\item A $G$ module is $\FF$-projective if and only if it is 
$\PP$-projective.
\item $\FF\cohom*G- \cong \PP\cohom*G-$
\end{enumerate}

\end{theorem}

\noindent{\sl Proof:} (i) It is obvious that any $\FF$-split sequence
is $\PP$-split, and the converse follows from a standard averaging
argument. Let $H$ be an arbitrary finite subgroup of $G$. Then $|H|
=\prod_{i=1,...,n} p_i^{a_i}$ where $p_i$ are distinct primes and $0
< a_i \in \Z$.  For each $i$, let $n_i$ be the index $n_i=[H:P_i]$. 
Now consider a $\PP$-split surjection $A\Epi{\pi} B$.  Let
$\sigma_i$ be a $P_i$-splitting of $\pi$, and define a map $s_i$ by
summing $\sigma_i$ over the cosets of $P_i$:
$$s_i(b) = \sum_{t\in H/P_i} t \sigma_i(t^{-1}b).$$ 
For each $P_i$ we obtain a map
$s_i: B \to A$, such that $\pi\circ s_i=n_i\times id_B$.  There exist 
$m_i\in\Z$ so that $\sum_i m_in_i=1$, and the map $s=\sum_i m_is_i$ 
is the required $H$-splitting.

\medskip\noindent (ii) It is obvious that a $\PP$-projective module is 
$\FF$-projective.  Now let $P$ be $\FF$-projective.  We may take a 
$\PP$-split surjection $M\epi P$ with $M$ a $\PP$-projective.  By~(i) 
this surjection is $\FF$-split, and hence split.  Thus $P$ is a direct 
summand of a $\PP$-projective and so is $\PP$-projective.  

\medskip\noindent (iii) now follows directly from (i) and (ii). \qed

\begin{proposition}\label{ffp0}
A group $G$ is of type $\FF\FP_0$ if and only if $G$
has only finitely many conjugacy classes of subgroups of prime power
order.  
\end{proposition}

\noindent {\sl Proof:}  Suppose that $G$ has only finitely many 
conjugacy classes of subgroups in $\PP$.  Let $I$ be a set of 
representatives for the conjugacy classes of $\PP$-subgroups and set
$$\Delta_0 = \bigsqcup_{P\in I} G/P .$$ This $G$-set satisfies
condition ($*$) for $\PP$ and therefore the surjection $\Z\Delta_0
\epi \Z$ is $\FF$-split and also $\Z\Delta_0$ is finitely generated.

To prove the converse we consider an arbitrary $\FF$-split surjection
$P_0 \epi \Z$ with $P_0$ a finitely generated $\FF$-projective. As in
\cite[6.1]{N1} we can show that $P_0$ is a direct summand of a module
$\bigoplus_{\delta\in \D_f} \ind{G_\delta}{G}{P_\delta}$, where $\D_f$
is a finite $G$-set, the $G_\delta$ are finite groups and $P_\delta$
are finitely generated $G_\delta$-modules. Therefore we might assume
from now on that $P_0$ is of the above form. Since there is a $G$-map
$\D_f \to \D$, where $\D$ satisfies condition $(\ast)$ the $\FF$-split
surjection $P_0 \Epi{\varepsilon} \Z$ is also $\D_f$-split
\cite{N1}. Consider now the following commutative diagram:

$$\xymatrix {P_0 \ar@{->>}[r]^{\varepsilon} \ar@/^/[d]^{\beta} &  \Z \ar@{=}[d] \\
                \Z\Delta_f  \ar@{->>}[r]^{\varepsilon_f} \ar@/^/[u]^{\alpha} & \Z }$$

That we can find such an $\alpha$ follows from the fact that
$\varepsilon$ is $\D_f$-split, and $\beta$ exists since $P_0$ is
$\D_f$-projective being a direct sum of induced modules, induced from
$G_\delta$, ($\delta\in \D_f$) to $G$.

As a next step we'll show that $\varepsilon_f$ is $\FF$-split. Take an
arbitrary finite subgroup $H$ of $G$ and show that $\varepsilon_f$
splits when restricted to $H$. Since $\varepsilon$ is split by $s$,
say, when restricted to $H$ we can define the required splitting by
$\beta \circ s$.

Now let $P$ be an arbitrary $p$-subgroup of $G$. Since the module $\Z [G/P]$
is $\FF$-projective, there exists a $G$-map $\varphi$, such that the following
diagram commutes:

$$\xymatrix{\Z\Delta_f  \ar@{->>}[r]^{\varepsilon_f} & \Z \ar@{=}[d] \\
                \Z [G/P] \ar@{->>}[r] \ar[u]^{\varphi}  & \Z}$$

The image $\varphi(P)$ of the identity coset $P$ is a point of
$\Z\Delta$ fixed by the action of $P$.  If $H$ is any group and 
$\Z\Omega$ is any permutation module, then the $H$-fixed points 
are generated by the orbit sums $H.\omega$.  Hence $P$ must 
stabilize some point of $\Delta_f$, since otherwise we would have 
that $p$ divides $\epsilon_f\varphi(P)=\epsilon\alpha(P)=1$, a 
contradiction.  It follows that $P$ is a subgroup of $G_\delta$ 
for some $\delta\in\Delta_f$.  
\qed


Note that being of type $\FF\FP_0$ does not imply that there are
finitely many conjugacy classes of finite subgroups.  In fact, the
authors have examples with infinitely many conjugacy classes of finite
subgroups, see \cite{vfg}.  Nevertheless this gives rise to the
following conjecture:

\begin{conjecture}
A group $G$ is of type $\FF\FP_\infty$ if and only if $G$ is of type 
$\FP_\infty$ and has finitely many conjugacy classes of $p$-subgroups.
\end{conjecture}

It is shown in~\cite{N1} that any $G$ of type $\FF\FP_\infty$ is of
type $\FP_\infty$, which together with Proposition~\ref{ffp0} proves 
one implication in the above conjecture.

\end{document}